\documentclass[11pt]{article} 

\usepackage[utf8]{inputenc} 

\usepackage[left=1.2in,right=1.2in,top=1.2in,bottom=1.2in]{geometry}

\usepackage{graphicx} 

\usepackage[parfill]{parskip}
\usepackage{color}

\usepackage{booktabs} 
\usepackage{array} 
\usepackage{verbatim} 
\usepackage{subfig} 

\usepackage{fancyhdr} 
\title{Applied Algebra and Geometry: A SIAGA of Seven Pictures}
\author{Anna Seigal \\ University of California, Berkeley}

\usepackage{amsmath}
\usepackage{amssymb}
\usepackage{amsthm}
\usepackage{mathrsfs}
\usepackage{mathtools}
\usepackage{tikz}
\usepackage[hidelinks]{hyperref}

\def\R{{\mathbb R}}

\def\Z{{\mathbb Z}}
\def \P{{\mathbb P}}

\begin{document}

\maketitle

\begin{abstract}

The cover of the SIAM Journal on Applied Algebra and Geometry shows seven
pictures. We describe these pictures and discuss the topics they represent.

\emph{About the Author:} Anna Seigal is a graduate student at UC Berkeley working in applied algebra. Her interests lie in tensors and applications to biological systems. She has a particular penchant for writing about mathematics in terms of pictures and also blogs on this subject at \url{https://picturethismaths.wordpress.com/}.

\end{abstract}


\section{Polynomial Optimization}

\begin{figure}[h]
\includegraphics[width=8cm]{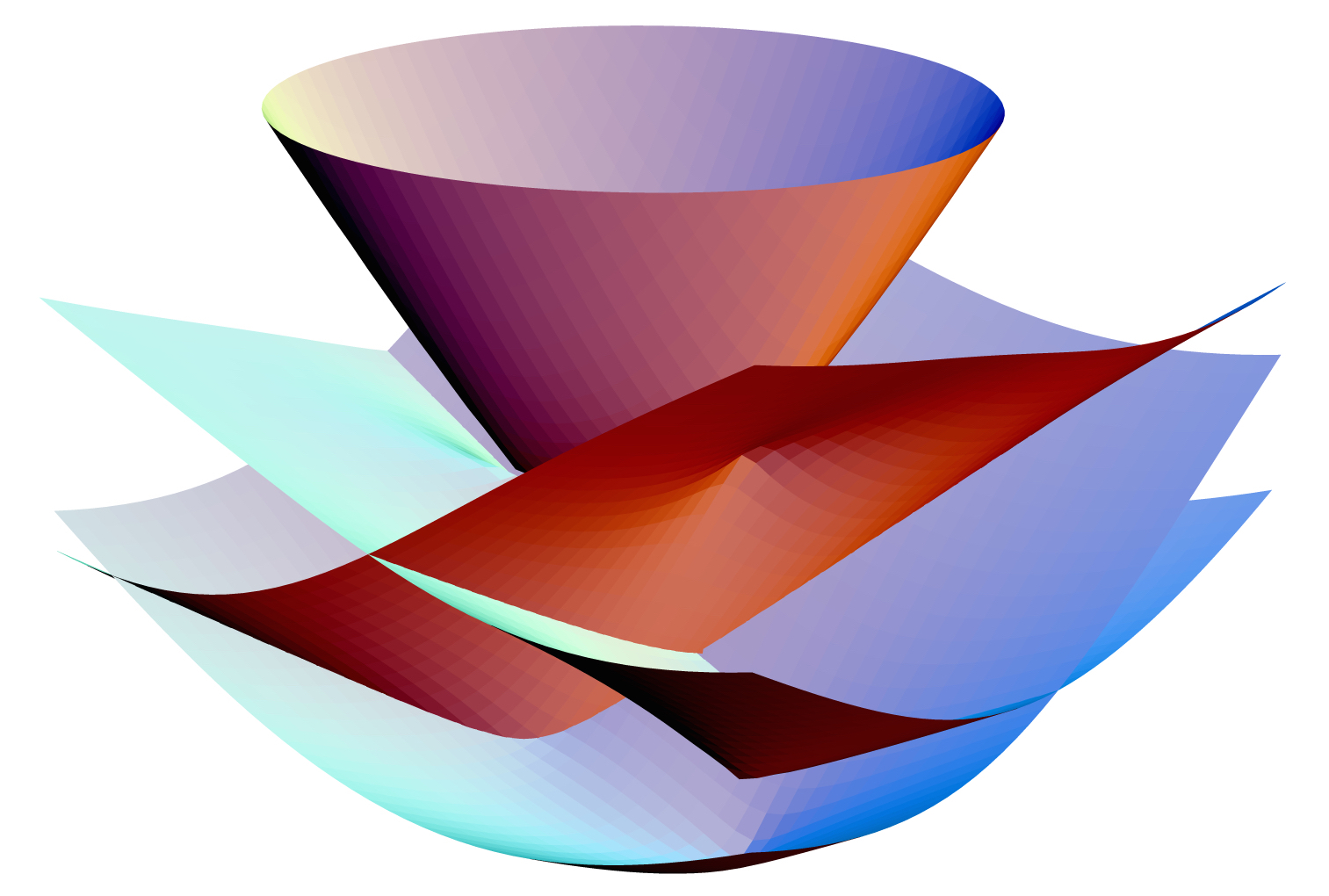}
\centering
\end{figure}

\subsubsection*{The Context}

Optimization is a central pillar of applied mathematics, with applications in fields from biology to engineering to finance. Polynomial optimization involves maximizing or minimizing a polynomial function subject to constraints given by polynomial equations. The feasible region, or the set of points in space that satisfy the required constraints, is our geometric object of interest. When considering particular structure on the constraints of an optimization problem, we have corresponding information on the kind of shape the feasible region can take. The boundary of the feasible region is of particular interest, since this is where the optimizing solution will often be found.

The geometry of a feasible region elucidates important aspects of the optimization problem. For example, it may indicate the type of algorithm best suited to maximize a function on that shape. In return, the family of geometrical shapes naturally associated with an optimization problem connects it to the expertise of other areas of mathematics.

Semi-definite optimization is a generalization of linear optimization. Here, we extend from constraints given by positivity conditions on the entries of a vector to the notion of positivity for a matrix. The constraints take the form of positive semi-definiteness of real symmetric matrices. Feasible regions take the form of spectrahedra, which result from intersecting the space of positive definite matrices with a linear space. Semi-definite optimization problems are often a valuable stepping-stone in understanding more complicated problems, and can be used as a \emph{relaxations} of the harder problems. For more on semi-definite optimization and its connection to algebraic geometry, see \cite{BPT}. 

\subsubsection*{The Picture}

A naturally-occurring feasible region with a nonlinear boundary is the circle
\[ (x- u_1)^2 + (y - v_1)^2 = d^2 , \]
the collection of points a fixed distance away from the center.

Perhaps instead we want to be a fixed distance away from two points - close to both the train station and the ferry terminal, for example. The collection of points whose sum of distances from two points is a constant describes an ellipse
\[ \left\{ (x, y) \in \R^2 : \sum_{i = 1}^2 \sqrt{ (x - u_i)^2 + (y - v_i)^2 } = d \right\} \text{.} \]

For problems involving, say, three factories, five train stations and a ferry terminal, we want to generalize this notion to the so-called $k$-ellipse. This is the set of points whose sum of distances from $k$ given points is equal to some constant.

We fix the locations of $k$ focal points $(u_i, v_i)$ and consider the distance $d$ to be an unknown. The picture shows the surface
\[ \left\{ (x, y,d) \in \R^3 : \sum_{i = 1}^k \sqrt{ (x - u_i)^2 + (y - v_i)^2 } = d \right\} \]
of points for which the sum of the distances from $(x,y)$ to the focal points is $d$. A polynomial equation can describe this constraint, which is the nonlinear boundary of the feasible region in a particular semi-definite optimization problem. 

The central convex part of the picture holds the solution to this minimization problem. Its lowest point is the Fermat-Weber point, which is often sought. If the distance of interest, $d$, is fixed, then we take a horizontal slice through the picture and optimize on this slice. The external components are also algebraic solutions to the polynomial constraints.

This picture first appeared in \cite{NPS}, and this version of the image is due to Cynthia Vinzant. It appeared as cover art for the May 2014 issue of the American Mathematical Society's Notices.

\section{Robotics}

\begin{figure}[h]
\includegraphics[width=8cm]{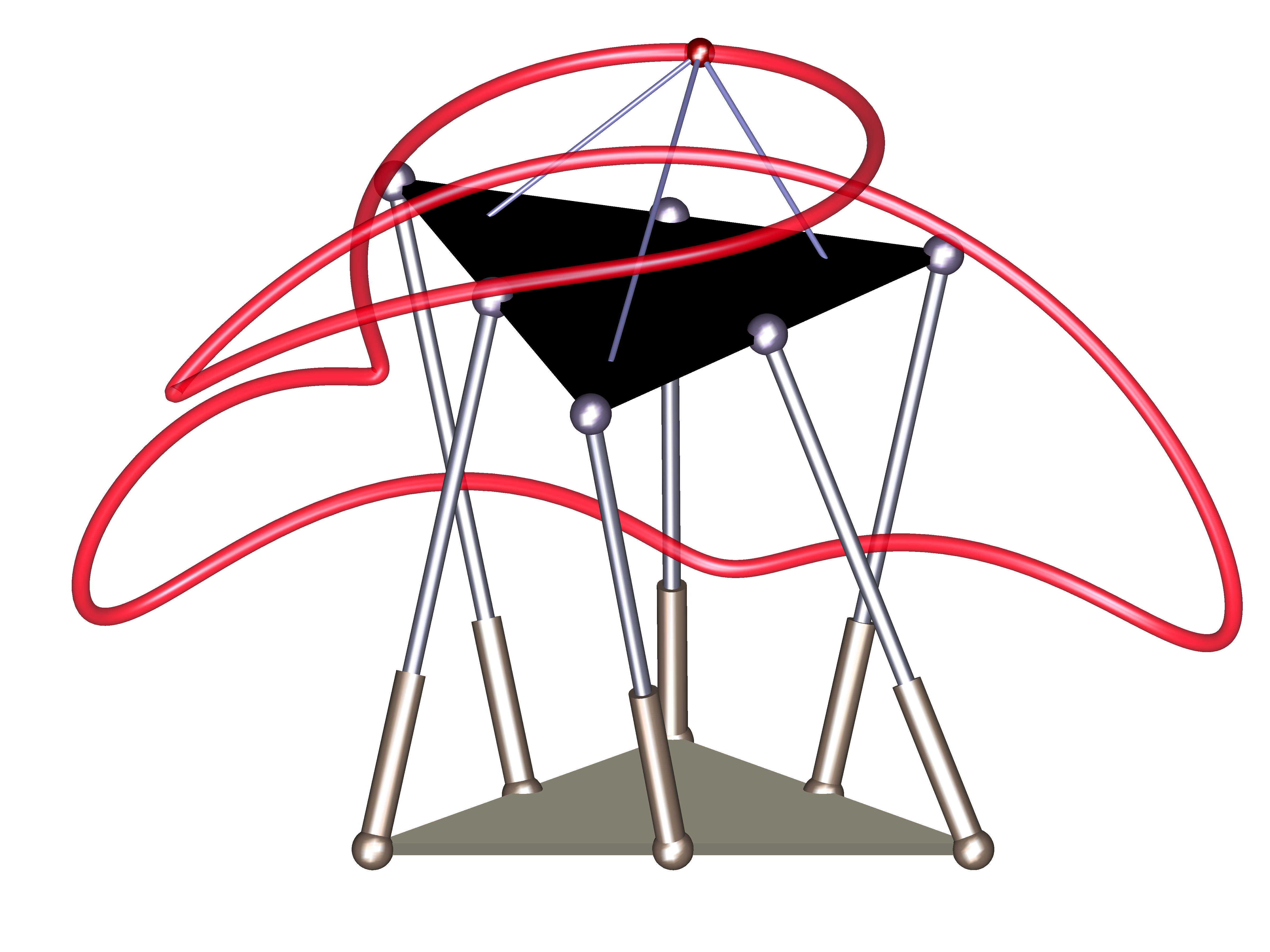}
\centering
\end{figure}

\subsubsection*{The Context}

Robotics---the design of mechanical machines to perform complex tasks---is a burgeoning field of study. The range and precision of robots' motion are limited by the mechanics of their constituent parts, which traditionally consist of rigid pieces connected by joints. The interplay between robotics, algebra, and geometry arises naturally from the kinematics of these pieces working together. In fact, numerical algebraic geometry largely arose from these applications. For an introduction, see Chapter 6 of \cite{CLO}.

A rigid part floating unconstrained in three-dimensional space has six degrees of freedom, but the joints of a mechanism restrict its motion. For most joints used in robotics, polynomial equations can describe these restrictions. With multiple pieces working together, one task is to find the location of a selected terminating part. For example, for fixed locations and angles of your shoulder, elbow, wrist and carpometacarpal joint, what is the location of your thumb? This is known as \emph{forward kinematics}. Furthermore, given a desired location of your thumb, what are the possible angles of your arm and hand that would achieve that location. This is \emph{inverse kinematics}.

The numerical algebraic geometry tool of \emph{homotopy continuation} can help answer these questions. In homotopy continuation, we first solve an easier but related set of polynomial equations. Then we deform the easier system to the more challenging problem of interest via a homotopy map. We use the solutions of the easier problem to obtain those of the harder problem by numerically tracking the paths of the original solutions as they deform under the homotopy. For more, see \cite{BHSW}.

\subsubsection*{The Figure}

The image depicts a special combination of rigid pieces and joints called a \emph{Griffis-Duffy Type I Platform}. It consists of two equilateral triangles, one fixed at the base and the other held above it by six rigid legs. Each leg connects a vertex of one triangle to a midpoint of an edge on the other. Although the lengths of each leg are fixed, the angles at each joint are free to move.

The geometry of this problem yields a system of polynomial equations that describes its kinematics: if we fix the point shown on top of the upper triangle, the collection of positions it can reach is the red curve, which is an algebraic set of degree 40.

The picture was created by Charles Wampler of General Motors and Douglas Arnold
of the University of Minnesota. It appeared on the poster for the IMA Thematic Year on Applications of Algebraic Geometry in 2006-07, in which significant progress was made connecting the use of algebraic geometry tools to industrial and applied mathematics.

\section{Polyhedral Geometry}

\begin{figure}[h]
\includegraphics[width=8cm]{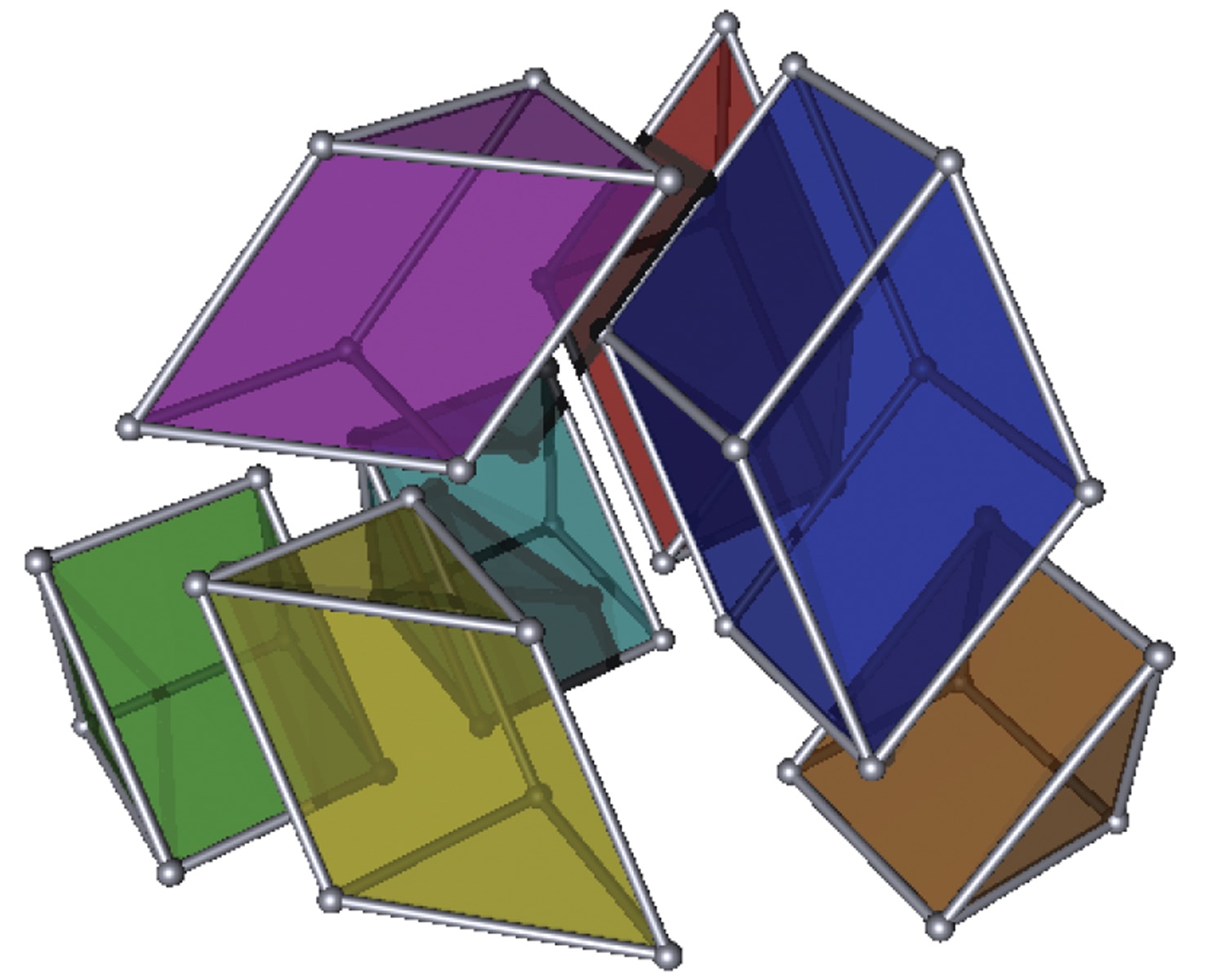}
\centering
\end{figure}

\subsubsection*{The Context}

A polyhedron is a combinatorial object that crops up in many places. For example, it is the shape of the feasible region for a linear optimization problem. It is a convex shape in $d$-dimensional space $\R^d$ described by an intersection of finitely many closed half-spaces
\[ P = \{ \mathbf{x} \in \R^n : A \mathbf{x} \leq b \} \]
where $A$ is a $d \times n$ matrix describing the angles of the half spaces, and $b \in \R^d$ encodes their translational information.

A compatible collection of polyhedra is called a polyhedral complex. It is possible to associate a polyhedral complex to an algebraic variety, allowing us to use combinatorial tools to understand the variety. A typical way to do this is via the methods of toric geometry. This approach has been applied in areas such as phylogenetics, integer programming, economics, biochemical reaction networks and computer vision (from where this particular picture arose).

Tropical geometry gives one method, called tropicalization, for getting a polyhedral complex from an algebraic variety. The new object gives us useful information. For example, the dimension of the original variety equals that of the new polyhedral complex and the latter is much easier to compute.

Polyhedral geometry also arises at the interface with the life sciences. For instance, Gheorghe 
Craciun recently announced a proof of the Global Attractor Conjecture for toric dynamical systems. 
This result is important for systems biology, and polyhedral fans play a key role in the proof.

\subsubsection*{The Picture}

This picture describes one example of using polyhedral tools to understand an algebraic variety, this shedding light on the application from which the variety arose.

In the field of \emph{computer vision}, and in the real world, `taking a photo' is a map from the three-dimensional world to a two-dimensional world. As any good photographer knows, the resulting features of the photo depends heavily on the angle and location of the camera.

We photograph three-dimensional projective space $\P^3$. Each camera, $A$, is a $3 \times 4$ matrix which determines a map $\mathbf{x} \mapsto A\mathbf{x}$ to two-dimensional projective space $\P^2$. This map tells us where each point of the original world ends up in the photograph.

More information can be gained by considering multiple cameras $( A_1 , A_2, A_3)$ at different locations. Then we have a map from the real world to three photographs:
\[ \phi : \P^3 \dashrightarrow (\P^2)^3 \]
\[ \hspace{15.5ex} \mathbf{x} \longmapsto (A_1 \mathbf{x}, A_2 \mathbf{x}, A_3 \mathbf{x} ) \]

The closure of the image of this map is an irreducible variety. For example, if the $A_i$ are the coordinate projections, the variety in $(\P^2)^3$ is cut-out by the Gröbner basis
\[ \{ z_0 y_2 - x_0 z_2 , z_1 x_2 - x_1 z_2, z_0 y_1 - y_0 z_1, x_0 y_1 x_2 - y_0 x_1 y_2 \} \text{ .} \]
When we take the initial monomials of these generators, their zero-set decomposes into seven pieces: one copy of $\P^1 \times \P^1 \times \P^1$, and six copies of $\P^1 \times \P^2$.

Our picture shows the three-dimensional shape we get from this zero-set when we identify each projective space $\P^i$ with the $i$-simplex. For example, $\P^2$ corresponds to the two-dimensional simplex $\Delta_2$ --- the triangle --- under the map:
\[ \P^2 \ni (x_0 : x_1 : x_2) \longleftrightarrow \frac{1}{x_0 + x_1 + x_2} ( x_0, x_1, x_2 ) \in \Delta_2 \]
and we identify each copy of $\P^1$ with the one-dimensional simplex $\Delta_1$ --- the unit-length line. 

Our zero-set is represented by a collection of polytopes which are faces of $(\Delta_2)^3$. The $\P^1 \times \P^1 \times \P^1$ corresponds to $\Delta_1 \times \Delta_1 \times \Delta_1$. This is the dark blue cube. Each of the six pieces $\P^2 \times \P^1$ correspond to $\Delta_2 \times \Delta_1$, a triangle cross a line segment. This gives the six triangular prisms in the picture. Each piece has been separated a little to make it easier to see, but the close-by parallel faces show how the different parts meet. Meeting at a triangle $\Delta_2$ means the projective spaces meet at a copy of $\P^2$. If the shared facet is a square $\Delta_1 \times \Delta_1$, the projective spaces meet at a copy of $\P^1 \times \P^1$. The original picture, and other nice ones, can be found in \cite{AST}. It was made using Michael Joswig's software `Polymake' \cite{GJ}.

\section{Topology of Data}

\begin{figure}[h]
\includegraphics[width=8cm]{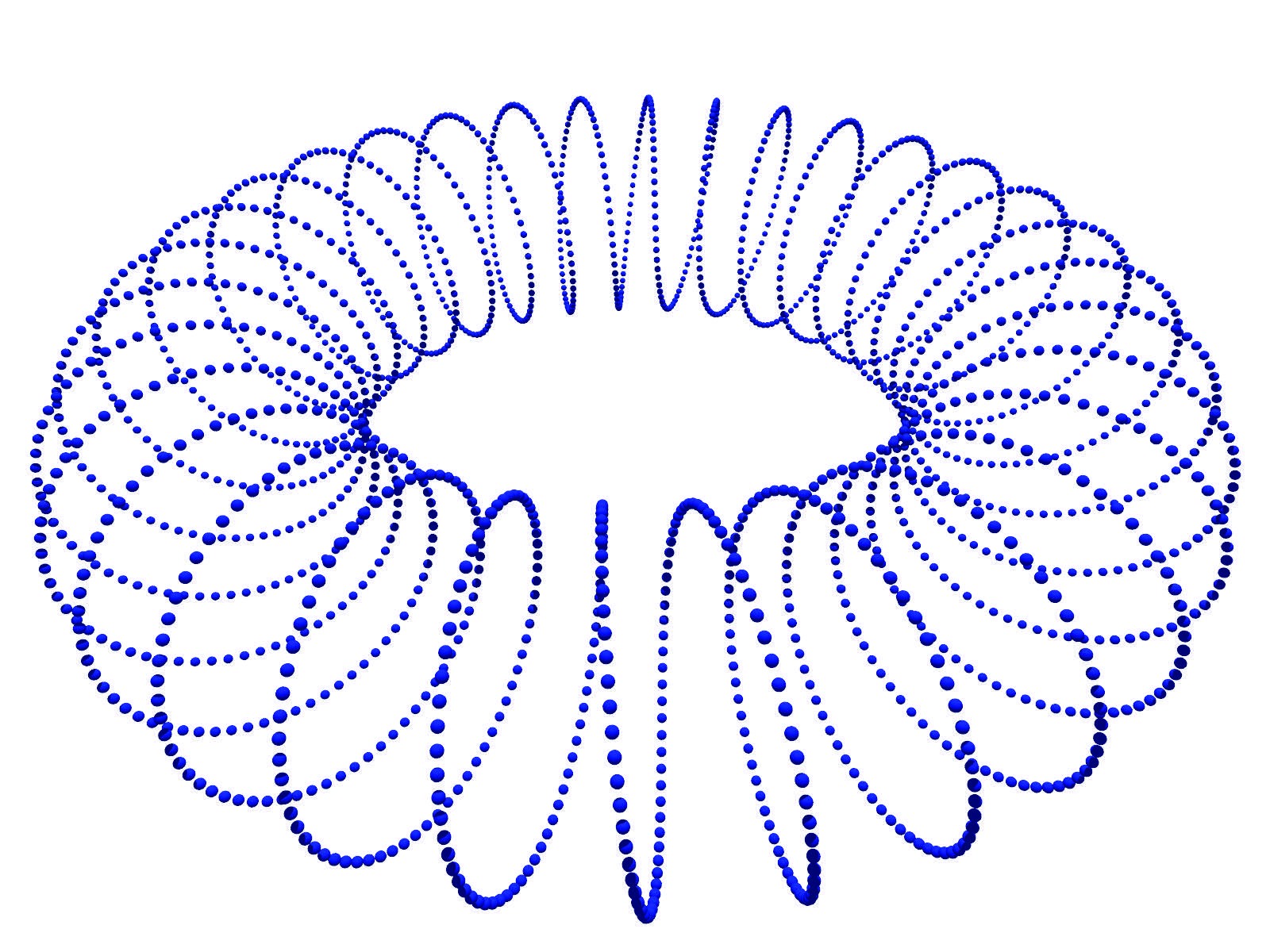}
\centering
\end{figure}

\subsubsection*{The Context}

Topology offers a set of tools that can be used to understand the shape of data. The techniques detect intrinsic geometric structures that are robust to many common sources of error including noise and arbitrary choice of metric. For an introduction, see \cite{Gh} and \cite{Ca}.

Say we have noisy data points coming from some unknown space $X$ which we believe possesses an interesting shape. We are interested in using the data to capture the \emph{topological invariants} of the unknown space. These are its holes of different dimensions, unchanged by continuous squeezing and stretching.

The holes of different dimensions are the homology groups of the space $X$. They are denoted by $H_k(X)$ for $k$ some non-negative integer. The zeroth homology group tells us the number of zero-dimensional holes or, more intuitively, the connectedness of the space. For a space $X$ with $n$ connected components, it is
\[ H_0(X) = \Z^n, \]
the free abelian group with $n$ generators. One-dimensional holes are counted by $H_1(X)$. For example, a circle $X = S^1$ has a single one-dimensional hole, so $H_1(S^1) = \Z$.

The connectedness properties of sampled data tell us a lot about the underlying space from which they are sampled. In some situations, such as for structural biological information, it is indispensable to know the structure of the holes too. These features are unchanged no matter which metric we use, or which space we embed the points into. The higher homology groups $H_k(X)$ for $k \geq 2$ similarly give us such summarizing features.

But there's a problem: sampling $N$ points from a space gives us a collection of zero-dimensional pieces, which --- unless two points land in exactly the same place --- are all unconnected. Let us call this data space $D_0$. The space $D_0$ has homology groups
\[ H_k(D_0) = \begin{cases} \Z^N & k = 0 \\ 0 & \text{otherwise.} \end{cases} \]

It is usually the case that many points are very close together, and ought to be considered to come from the same connected component. To measure this we use \emph{persistent homology}. We take balls of increasing size centered at the original data points, and measure the homology groups of the space consisting of the union of these balls. We call this space $D_\epsilon$, where $\epsilon$ is the radius of the balls. The important structural features are those that \emph{persist} for large ranges of values of $\epsilon$.




\subsubsection*{The Picture}

This picture shows data points sampled from a torus, which we imagine to live in three-dimensional space. It was made by Dmitriy Morozov, who works at the Lawrence Berkeley National Lab. He applies topological methods in cosmology, climate modeling and material science.

The sampled points in the picture lie on the torus, and furthermore in a more specialized slinky-shaped zone of the torus. This is an important feature of the shape which topological methods will capture.


The original data consists of 5000 points, and our persistent homology approach involves taking three-dimensional balls $B_\epsilon(d_i)$ of radius $\epsilon$ centered at each data point $d_i$. When the radius $\epsilon$ is very extremely small, none of the balls will be connected, and the shape of our data is indistinguishable from any other collection of 5000 points in space.

Before long, the radius will exceed half the distance to all the points' nearest neighbors. The 5000 balls join together to form a curled up circular piece of string. Topological invariants do not notice the curling, so topologically the shape obtained is a thickened circle with a one-dimensional hole $H_1(D_{\epsilon}) = \Z$. When the radius is large enough for the adjacent curls of the slinky to meet, but not to read the opposite side of each curl, we get a hollow torus with $H_1(D_{\epsilon}) = \Z^2$ and $H_2(D_\epsilon) = \Z$. Finally, the opposite sides of each curl of the slinky will meet, and they will meet up with the slinky-curls on the opposite side of the torus. Our shape then becomes a three-dimensional shape with no holes, and $H_1(D_R) = 0$.

In this example, the data points can be visualized and we are able to confirm that our intuition for the important structure of the shape agrees with the homological computations. For higher-dimensional examples it is these persistent features that will guide our understanding of the shape of the data.

\section{Geometric Modeling}

\begin{figure}[h]
\includegraphics[width=8cm]{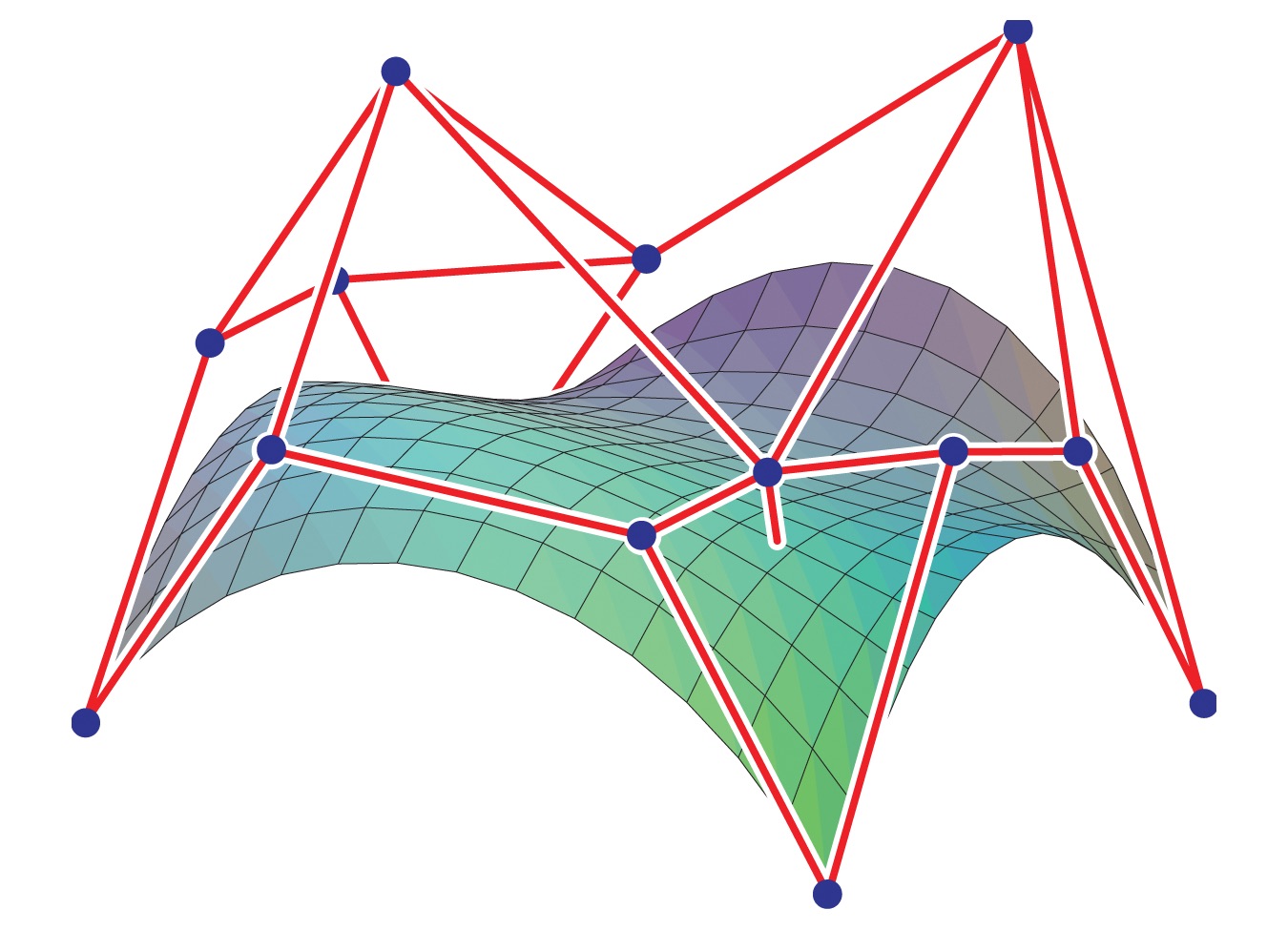}
\centering
\end{figure}
\subsubsection*{The Context}

Geometric Modeling is an area of applied mathematics in which piecewise polynomials are used to build computer models for depicting and describing shapes in space. 

One tool that is used for such modeling is a parametric curve called a \emph{B\'{e}zier curve}. They are named after the French engineer Pierre B\'{e}zier who worked in the automotive industry (like Charles Wampler from Section 2).

A B\'{e}zier curve models smooth motion through time or space. Each one is defined by a number of \emph{control points}: points which specify its shape and location. This make them easy to manipulate on a computer interface: changing the location of the control points causes a reliable change in the curve.

A collection of $d+1$ control points $P_0, \ldots, P_d$ defines a B\'{e}zier curve of degree $d$. The simplest example is when the degree is one, and we have two control points $P_0$ and $P_1$. In this case, the B\'{e}zier curve is the line that connects the two points
\[ B(t) = (1-t) P_0 + t P_1 \]
for $t$ between $0$ and $1$. A degree two example is given by 
\[ B(t) = (1-t)^2 P_0 + 2t(1-t) P_1 + t^2 P_2 \]
and in general we have
\[ B(t) = \sum_{i = 0}^d {{d \choose i}} {(1-t)}^{d-i} t^i P_i \text{ .} \]

A B\'{e}zier curve has a `control polygon' associated to its control points, which is found by taking the line segments connecting adjacent control points. The convex hull of this control polygon contains the curve. The control polygon has many other useful properties for example in approximation of the curve.

\subsubsection*{The Picture}

This picture shows a generalization of the B\'{e}zier curves described above to a two-dimensional B\'{e}zier surface. It is from \cite{GSZ}.

It is very useful for applications to have a nice way to make smooth two-dimensional surfaces. For example they have been used in the design for parts of a car.

A B\'{e}zier surface is defined in terms of a collection of control points in three-dimensional space
\[ \{ P_{0,0}, \ldots, P_{d_1, d_2} \} \]
which now are indexed by two indices rather than one. It is given parametrically by
\[ B(t_1, t_2) = \sum_{i_1 = 0}^{d_1} \sum_{i_2 = 0}^{d_2} \left( {{d_1 \choose i_1}} {(1-t_1)}^{d_1-i_1} t_1^{i_1} \right) \left( {{d_2 \choose i_2}} {(1-t_2)}^{d_2-i_2} t_2^{i_2} \right)  P_{i_1, i_2} \text{ .} \]
A list of $(d_1 + 1)(d_2 + 1)$ control points gives a surface of degree $d_1 d_2$ via this map.

The control points are shown in blue. The B\'{e}zier surface is shown below them in green. We now have a two-dimensional analogue of the control polygon, below whose convex hull the surface sits. This is shown by red lines connecting the blue points. This polyhedral structure connects this topic to the one in Section 3.

Applications often demand the investigation of further properties of B\'{e}zier curves and surfaces, such as how they intersect with one another. One step in the process is: given a parametric description of a surface, obtain an implicit description of it. That is, find the relations amongst the coordinates that are satisfied for all points on the surface. Here, computational algebraic geometry tools are very useful.

\section{Tensors}

\begin{figure}[h]
\includegraphics[width=8cm]{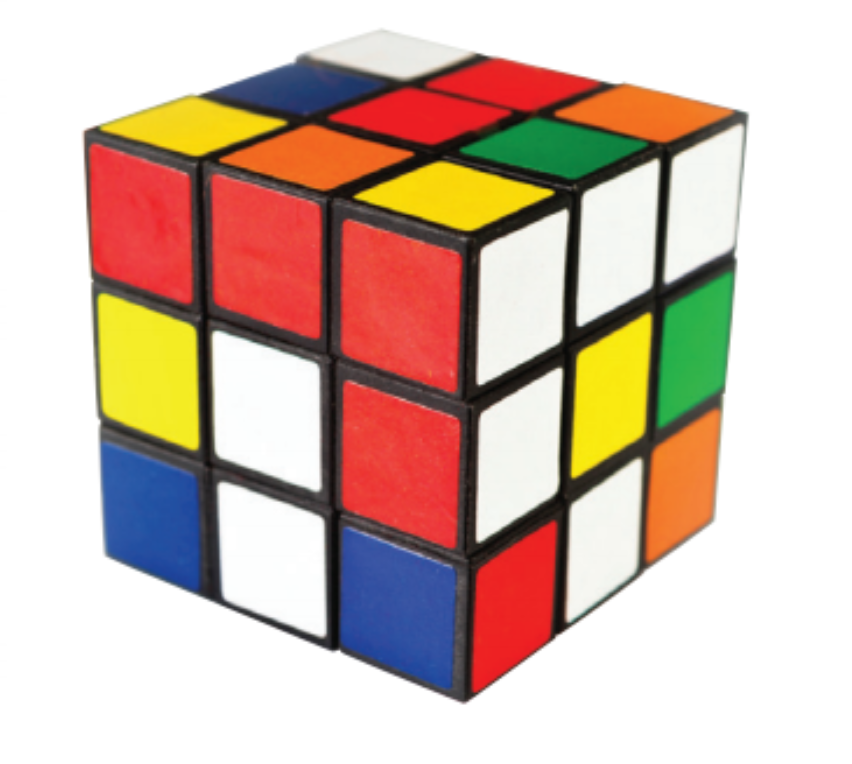}
\centering
\end{figure}

\subsubsection*{The Context}

Tensors are the higher-dimensional analogues of matrices. They are data arrays with three or more dimensions, and are represented by an array of size $n_1 \times \cdots \times n_d$, where $n_k$ is the number of `rows' in the $k$th direction of the array. The entries of the tensor $A$ are denoted by $A_{i_1 \ldots i_d}$ where $i_k \in \{ 1, \ldots, n_k \}$ tells you which row in the $k$th direction you are looking at. Just as for a matrix, the entries of a tensor are elements in some field, for example real or complex numbers.

Tensors occur naturally when it makes sense to organize data by more than two indices. For example, if we have a function that depends on three or more discretized inputs $f(x,y,z)$ where $x \in \{ x_1, \ldots, x_{n_1} \}$, $y \in \{ y_1, \ldots, y_{n_2} \}$ and $z \in \{ z_1, \ldots, z_{n_3} \}$, then we can organize the values $A_{ijk} = f(x_i,y_j,z_k)$ into a tensor of size $n_1 \times n_2 \times n_3$. Tensors are increasingly widely used in many applications, especially signal processing, where the uniqueness of a tensor's decomposition allows the different signals comprising a mixture to be found. They have also been used in machine learning, genomics, geometric complexity theory and statistics.

Our data analysis techniques are currently limited to a matrix-centric perspective. To overcome this, there has been tremendous effort to extend the well-understood properties of matrices to the higher-dimensional world of tensors. A greater understanding of tensors paves the way for very exciting new developments that can cater to the natural structure of tensor-based data, for example in experimental design or confounding factor analysis. This analysis and understanding uses interesting and complicated geometry.

One requirement for computability of a tensor is to have a good low rank approximation. Tensors of size $n_1 \times \cdots \times n_d$ have $n_1 \ldots n_d$ entries and, for applications, this quickly becomes unreasonably large. Matrices are analyzable via their singular value decomposition, and the best low rank approximation is obtainable directly from this by truncating at the $r$th largest singular value. For tensors we can also define useful related notions such as eigenvectors, singular vectors, and the higher order singular value decomposition.
 
\subsubsection*{The Picture}

As well as being a picture of the well-known Rubik's cube, this picture describes a cartoon of a tensor of size $3 \times 3 \times 3$.  Such a tensor consists of 27 values.

To understand the structure contained in a tensor, we use its natural symmetry group to find a presentation of it that is simple and structurally transparent. This motivation also underlies the Rubik's puzzle although the symmetries can be quite different: a change of basis transformation for the tensor case, and a permutation of pieces in the case of the puzzle.

Despite being small, a $3 \times 3 \times 3$ tensor has interesting geometry. It is known that a generic tensor of size $3 \times 3 \times 3$ has seven eigenvectors in $\P^2$. We show in \cite{ASS} that any configuration of seven eigenvectors can arise, provided no six of the seven points lie on a conic.

\section{Visualization of Algebraic Varieties}

\begin{figure}[h]
\includegraphics[width=8cm]{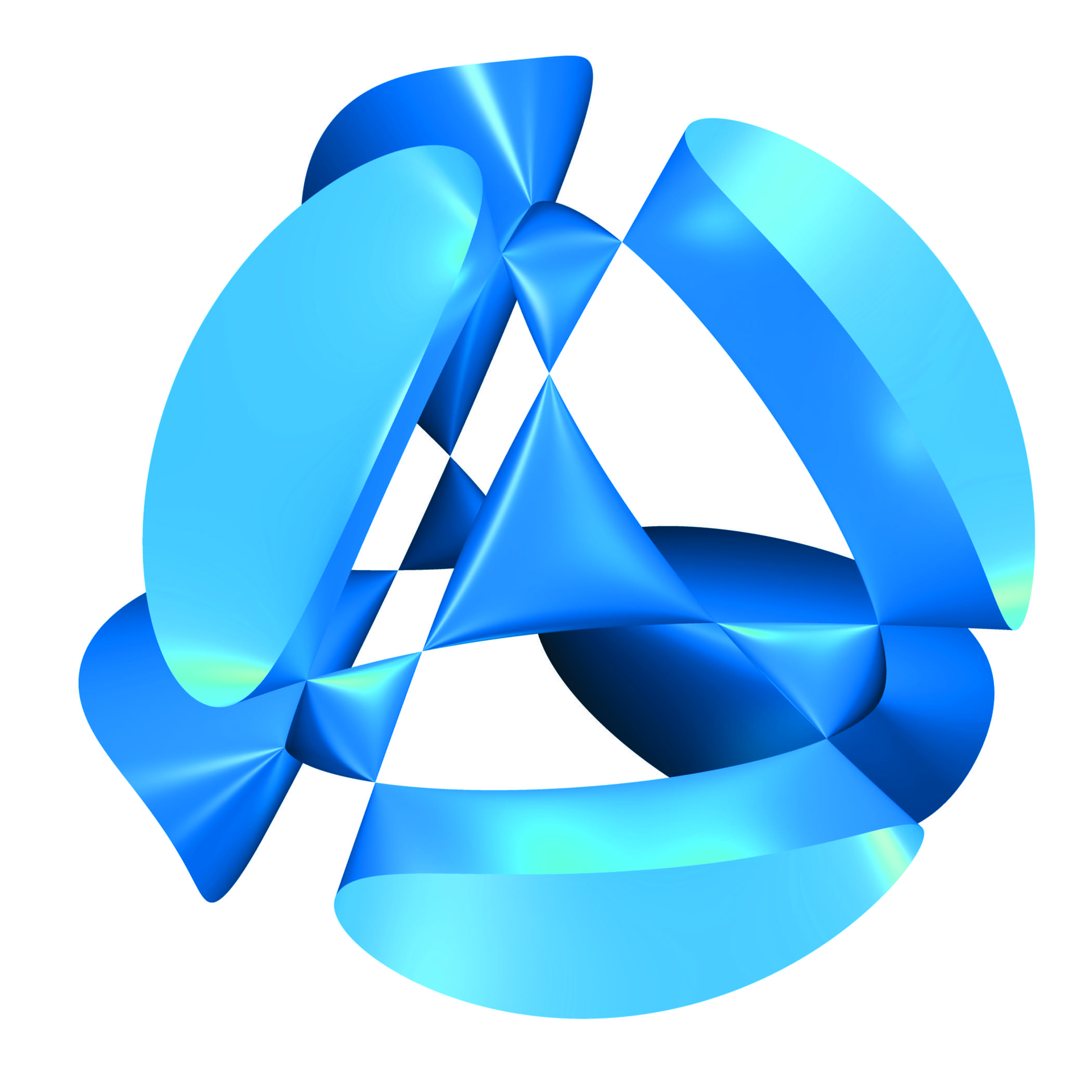}
\centering
\end{figure}

\vspace{-6ex}

\subsubsection*{The Context}

There is a vast mathematical toolbox of techniques that can be used to understand algebraic varieties. We have encountered some so far in this document, for example polyhedral geometry. It is great when we are actually able to draw the algebraic variety in question, using visualization software. When possible, this facilitates the most direct of observations to be made.

Although it poses an obvious restriction on the number of dimensions we can work in, even visualizing particular slices through our variety of interest is structurally revealing. Large polynomials with many terms can be very hard to get a handle on, and it makes sense to use modern-day computer tools to convert these equations into helpful pictures.

\subsubsection*{The Picture}

This picture shows a \emph{Kummer Surface}. It was made by Oliver Labs using the visualization software `Surfex'. Many beautiful pictures have been created in this way: for more, see the picture galleries from the `Imaginary: Open Mathematics' website. 

It is an example of an irreducible surface in three-dimensional space of degree four. In general, these have at most 16 singular points. Kummer surfaces are those that attain this upper bound. The 16 singular points represent the 2-torsion points on the Jacobian of the underlying genus 2 curve.

This picture also represents the problem-solving areas of coding theory and cryptography, in which there can be found a broad range of applied algebra and geometry. The group law on an elliptic curve is fundamental for cryptography. Similarly, the group law on the Jacobian of hyperelliptic curves has been used for cryptographic purposes, see \cite{BCHL} and \cite{BSZ}, the former is by Kristin Lauter from Microsoft Research who is president of the Association for Women in Mathematics (AWM).

\bibliographystyle{amsalpha}

\end{document}